\documentclass[10pt] {amsart}

\usepackage{amssymb}

\usepackage{graphicx}

\oddsidemargin=\evensidemargin
\catcode`\@=11
\long\def\@makefntext#1{
\protect\noindent \hbox to 3.2pt {\hskip-.9pt
$^{{\eightrm\@thefnmark}}$\hfil}#1\hfill}		

\def\ps@myheadings{\let\@mkboth\@gobbletwo		
\def\@oddhead{\hbox{}
\rightmark\hfil\eightrm\thepage}
\def\@oddfoot{}\def\@evenhead{\eightrm\thepage\hfil
\leftmark\hbox{}}\def\@evenfoot{}
\def\sectionmark##1{}\def\subsectionmark##1{}}

\catcode`@=11		      		     
\def\ps@plain{\let\@mkboth\@gobbletwo
     \def\@oddhead{}\def\@oddfoot{\eightrm\hfil\thepage
     \hfil}\def\@evenhead{}\let\@evenfoot\@oddfoot}

\newtheorem{theorem}{Theorem}[section]

\newtheorem{proposition}[theorem]{Proposition}
\newtheorem{corollary}[theorem]{Corollary}
\newtheorem{definition}[theorem]{Definition}

\numberwithin{equation}{section}

\def\abstracts#1#2#3#4{{
	\centering{\begin{minipage}{4.5in}\footnotesize\baselineskip=10pt
	\centerline{ABSTRACT}
	\parindent=15pt #1\par
	\parindent=15pt #2\par
	\parindent=15pt #3\par
	\parindent=15pt #4\par
	\end{minipage}}\par}}

\newcommand{\textlineskip}{\baselineskip=13pt}
\newcommand{\smalllineskip}{\baselineskip=10pt}

\renewenvironment{thebibliography}[1]
	{\frenchspacing
	 \ninerm\baselineskip=11pt
	 \begin{list}{[\arabic{enumi}]}
	{\usecounter{enumi}\setlength{\parsep}{0pt}
	 \setlength{\leftmargin 19pt}{\rightmargin 0pt}   
	 \setlength{\itemsep}{0pt} \settowidth
	{\labelwidth}{[#1]}\sloppy}}{\end{list}}

\newcommand{\fcaption}[1]{
        \refstepcounter{figure}
        \setbox\@tempboxa = \hbox{\footnotesize Fig.~\thefigure. #1}
        \ifdim \wd\@tempboxa > 5in
           {\begin{center}
        \parbox{5in}{\footnotesize\smalllineskip Fig.~\thefigure. #1}
            \end{center}}
        \else
             {\begin{center}
             {\footnotesize Fig.~\thefigure. #1}
              \end{center}}
        \fi}
\def\runninghead#1#2{\pagestyle{myheadings}
\markboth{{\protect\footnotesize\it{\quad #1}}\hfill}
{\hfill{\protect\footnotesize\it{#2\quad}}}}

\font\ninerm=cmr9

\font\eightrm=cmr8

\newcommand{\TryPackage}[3]{\IfFileExists{#1.sty}{\usepackage{#1}#2}{#3}}
\TryPackage{mathrsfs}{\renewcommand{\mathcal}{\mathscr}}{%
     \TryPackage{eucal}{}{}}

\date{\today}
\begin{document}
\setlength{\textheight}{7.7truein}  

\runninghead{\quad The Jones polynomial and boundary slopes}
{The Jones polynomial and boundary slopes \quad}

\normalsize\textlineskip
\thispagestyle{empty}
\setcounter{page}{1}

\vspace*{0.88truein}

\centerline{\bf THE JONES POLYNOMIAL AND BOUNDARY SLOPES}
\baselineskip=13pt
\centerline{\bf OF ALTERNATING KNOTS}
\vspace*{0.37truein}
\centerline{\footnotesize CYNTHIA L. CURTIS}
\baselineskip=12pt
\centerline{\footnotesize\it Department of Mathematics \& Statistics}
\baselineskip=10pt
\centerline{\footnotesize\it The College of New Jersey}
\baselineskip=10pt
\centerline{\footnotesize\it Ewing, NJ}
\baselineskip=10pt
\centerline{\footnotesize\it 08628}
\baselineskip=10pt
\centerline{\footnotesize\it {\tt ccurtis@tcnj.edu}}

\vspace*{10pt}
\centerline{\footnotesize SAMUEL TAYLOR}
\baselineskip=12pt
\centerline{\footnotesize\it Department of Mathematics \& Statistics}
\baselineskip=10pt
\centerline{\footnotesize\it The College of New Jersey}
\baselineskip=10pt
\centerline{\footnotesize\it Ewing, NJ}
\baselineskip=10pt
\centerline{\footnotesize\it 08628 USA}
\baselineskip=10pt
\centerline{\footnotesize\it {\tt staylor@math.utexas.edu}}


\vspace*{0.21truein}
\abstracts{We show for an alternating knot the minimal boundary slope of an essential spanning surface is given by the signature plus twice the minimum degree of the Jones polynomial and the maximal boundary slope of an essential spanning surface is given by the signature plus twice the maximum degree of the Jones polynomial. For alternating Montesinos knots, these are the minimal and maximal boundary slopes.
}{}{}{}

\vspace*{0.225truein}
\section*{Correction added January, 2014}
After the publication of this paper, Joshua Howie pointed out to us that Theorem $2.8$ incorrectly states a theorem of Adams and Kindred appearing in \cite{AK}. See \cite{How} for relevant examples and further explanation. Because of this error, Theorem $3.1$ is false as stated. In particular, for an alternating knot $K$ our result applies only to boundary slopes coming from the basic layered surfaces associated to $K$. The correct version of our main theorem is the following:

\begin{theorem}[Corrected version of Theorem 3.1]
Let $K$ be an alternating knot. Then the largest boundary slope of a basic layered surface for $K$ is twice the maximum degree of the Jones polynomial plus $\sigma(K)$ and the smallest boundary slope of a basic layered surface for $K$ is twice the minimum degree of the Jones polynomial plus $\sigma (K)$, where $\sigma(K)$ denotes the signature of the knot $K$.
\end{theorem}

To see that this is the case, note the correct statement of Corollary 2.9 is as follows:

\begin{corollary}[Corrected version of Corollary 2.9]
If $K$ is an alternating knot with reduced alternating diagram $D$, then the two checkerboard surfaces for $K$ given by $D$ have the maximal and minimal boundary slopes among all basic layered surfaces.
\end{corollary}

The proof of the main theorem follows immediately using this version of the Corollary. Finally, by results of Ichihara and Mizushima in \cite{IM1,IM2}, summarized in Theorem $2.10$, for alternating Montesinos knots we still obtain our main corollary.

\begin{corollary}
Let $K$ be an alternating Montesinos knot. Then the largest boundary slope of $K$ is twice the maximum degree of the Jones polynomial plus $\sigma(K)$ and the smallest boundary slope of $K$ is twice the minimum degree of the Jones polynomial plus $\sigma (K)$, where $\sigma(K)$ denotes the signature of the knot $K$.
\end{corollary}

We would like to thank Joshua Howie for pointing out our error. Below appears our original paper as it previously appeared.

\section{Introduction}
Let $K$ be a knot in $S^3$ and let $N(K)$ be a tubular neighborhood of $K$. Let $M = S^3 - N^{\circ}(K)$. The collection of essential surfaces in $M$ has proved to be a useful tool in understanding both the knot $K$ and the collection of 3-manifolds resulting from Dehn surgeries on $K$.  In recent years, attention has focused in particular on the slopes of the boundary curves of the essential surfaces with boundary in the boundary of $M$.  These slopes have played a key role in attempts to understand the non-hyperbolic, or exceptional, Dehn surgeries on hyperbolic knots in particular.  The slopes also play an important role in understanding the $A$-polynomial of $K$.

Much is known about the set of boundary slopes of certain classes of knots. By Hatcher \cite{H}, the number of boundary slopes is finite for any knot $K$. In \cite{HT}, Hatcher and Thurston calculate all boundary slopes for $2$-bridge knots, and in \cite{HO} Hatcher and Oertel generalize this to compute the boundary slopes for all Montesinos knots. More recently, Mattman, Maybrun, and Robinson proved that for $2$-bridge knots the difference between the largest and smallest boundary slope is equal to twice the crossing number of the knot \cite{MMR}. Similarly, Ichihara and Mizushima proved in \cite{IM1} that this difference is bounded above by twice the crossing number for any prime Montesinos knot, with equality for alternating Montesinos knots. For these knots, Ichihara and Mizushima show in \cite{IM2} that the largest and smallest boundary slopes are twice the number of positive crossings in a reduced alternating diagram for the knot and the negative of twice the number of negative crossings in a reduced alternating diagram for the knot, respectively.

Most recently, in \cite{AK}, Adams and Kindred have classified the spanning surfaces for alternating knots.  It is an easy consequence of their work that the largest and smallest boundary slopes of essential spanning surfaces are attained by the checkerboard surfaces arising from a reduced alternating diagram of the knot; we establish this below in Corollary 2.9. It follows from Proposition 2.5 that the largest and smallest boundary slopes of essential spanning surfaces are twice the number of positive crossings in a reduced alternating diagram for the knot and the negative of twice the number of negative crossings in a reduced alternating diagram for the knot, respectively.

It is also well-known that the difference between the highest and lowest degrees of terms of the Jones polynomial is equal to the crossing number of the knot for an alternating knot. In this paper, we show that the difference between twice the highest degree of the Jones polynomial and the largest boundary slope of an essential spanning surface and similarly the difference between twice the lowest degree of the Jones polynomial and the smallest boundary slope of an essential spanning surface is an invariant for alternating knots, namely the signature of the knot. Thus, for alternating knots, the Jones polynomial detects the maximal and minimal boundary slopes of essential spanning surfaces. In particular, for alternating Montesinos knots, the Jones polynomial detects the maximal and minimal boundary slopes in light of \cite{IM2}.

We note that recent related results appear in \cite{O} and \cite{FKP}.

\section{Checkerboard surfaces and maximal and minimal boundary slopes}
In this section we establish terminology and results regarding essential surfaces in knot complements and their boundary slopes. We begin by establishing terminology regarding essential surfaces in $M$.

A surface $S \subset M$ is said to be \textit{incompressible} if for any disc $D \subset M$ with $D \cap S = \partial D$, there exists a disc  $D' \subset S$, with $\partial D' = \partial D$.  A surface $S$ is $\partial${\em - incompressible} if for each disc $D \subset M$ with $D \cap S = \partial_+D$ and $D \cap \partial M = \partial_-D$ there is a disc $D' \subset S$ with $\partial_+ D' =\partial_+ D$ and $\partial_- D'  \subset \partial S$. A surface $S\subset M$ is \textit{essential} if it is both incompressible and $\partial$-incompressible.

If $S$ is properly embedded in $M$, then $\partial S$ is a disjoint collection of copies of $S^1$, each embedded in $\partial M$. It follows that the boundary components of $S$ are parallel curves on the boundary torus $\partial M$. Let $\mu$ and $\lambda$ denote the meridian and longitude of $K$, respectively. Then if $S$ is an essential surface and $C$ is one of the parallel components of $\partial S$, we may write $[C] = m[\mu] + l[\lambda]$ in $\pi_1 (\partial M)$.
\begin{definition}
The fraction $\frac{m}{l}$ is said to be the {\em boundary slope} of $S$ in $M$.
\end{definition}

We now describe the essential surfaces which realize the maximal and minimal boundary slopes for alternating Montesinos knots and the maximal and minimal integral boundary slopes for arbitrary alternating knots. Let $K$ be any knot positioned to lie in a plane $P = S^2 - \infty$ except in a neighborhood of each crossing, and let $D$ be the diagram of $K$ given by the projection of $K$ into $P$. Note that the projection of $K$ in $P$ divides $P$ into a finite number of regions, with four regions meeting at each crossing. If one arbitrarily chosen region is colored black, we may proceed to alternately color the regions given by the projection black or white so that at each crossing there are two regions of each color, with like-colored regions meeting diagonally at the crossing but not along an edge. Finally, the black regions may be joined along half-twisted bands at each crossing to give a surface with boundary $K$.  Similarly the white regions may be joined by half-twisted bands to form a second surface with boundary $K$.
\begin{definition}
These two surfaces are the {\em checkerboard surfaces} for $K$ given by the diagram $D$.
\end{definition}

Now if $K$ is an alternating knot and if $D$ is a reduced alternating diagram for $K$, we have
\begin{theorem}(Menasco and Thistlethwaite \cite{MT}, Proposition 2.3)
The checkerboard surfaces for $K$ given by $D$ are both essential.
\end{theorem}

 Orient $K$. Let $S$ be a checkerboard surface for $K$ given by $D$.
 \begin{definition}
 For each crossing of $D$, we say that the crossing is {\em nonexceptional} with respect to $S$ if the orientation of the two strands of $K$ at the crossing are opposite, so that they give rise to an orientation of the half-twisted band which is a piece of $S$ in a neighborhood of the crossing. If the two strands of $K$ have a parallel orientation along the half-twisted band, the crossing is {\em exceptional} with respect to $S$.
\end{definition}

The boundary slopes of the checkerboard surfaces of an alternating knot given by a reduced alternating diagram of the knot are well-known, but we compute them here for completeness.

\begin{proposition}
If $D$ is a reduced alternating diagram of $K$ then the boundary slopes of the checkerboard surfaces are $2cr_+$ and $-2cr_-$, where $cr_+$ (respectively $cr_-$)  is the number of positive (respectively, negative) crossings in $D$.
\end{proposition}

\begin{proof}
Let $S$ be a checkerboard surface given by $D$. The boundary slope of $S$ is the slope of the curve $C = S \cap N(K)$. This is the linking number of $C$ and $K$ in $S^3$, where $C$ and $K$ are given parallel orientations. This linking number is easily computed by considering the local contribution to linking number in a neighborhood of each crossing in $D$. Whether there is a contribution to linking at a given crossing depends on both the sign of the crossing and whether the crossing is exceptional. The four possible cases are shown in Figures 1 and 2. We see that for a given crossing the local contribution to linking is twice the sign of the crossing if the crossing is exceptional with respect to $S$ and zero otherwise. Hence the boundary slope of $S$ is twice the signed sum of the exceptional crossings with respect to $S$.

\begin{figure}
  \includegraphics[scale=1]{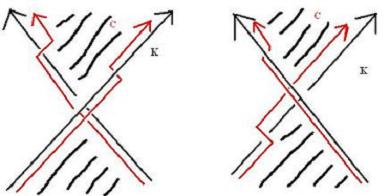}\\
  \caption{The local contribution to the crossing number is 2 or -2 for exceptional crossings, according to the sign of the crossing}\label{exceptcrossings}
\end{figure}
\begin{figure}
  \includegraphics[scale=1]{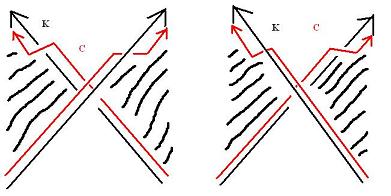}\\
  \caption{The local contribution to the crossing number is 0 for nonexceptional crossings}\label{nonexceptcrossings}
\end{figure}

Fix an orientation for $K$. Choose a crossing in $D$, and choose a strand of $K$ at the crossing.  Consider the next crossing, following the chosen strand of $K$ according to the orientation on $K$.  Since $D$ is alternating, it is easy to check that if the signs of the crossings agree, then the new crossing is exceptional if and only if the original crossing is exceptional.  However if the new crossing has the opposite sign from the original crossing, then the new crossing is exceptional if and only if the original crossing is nonexceptional. Proceeding similarly through all crossings of the knot, we see that all exceptional crossings of $K$ have the same sign, and all nonexceptional crossings have the opposite sign.  Hence the boundary slope of $S$ is $2cr_+$ if the exceptional crossings are positive or $-2cr_-$ if the exceptional crossings are negative.

Now let $S'$ be the other checkerboard surface given by $D$.  As for $S$, the boundary slope of $S'$ is either $2cr_+$ or $-2cr_-$. However note that at any crossing of $D$, if the crossing is exceptional for $S$ then it is nonexceptional for $S'$ and vice-versa by the definition of the checkerboard surfaces.  Thus, the positive crossings are the exceptional crossings for exactly one of the surfaces $S$ and $S'$, and the negative crossings are exceptional for the other surface. Hence once surface has boundary slope $2cr_+$ and the other has boundary slope $-2cr_-$.
\end{proof}

We remark that if $D'$ is another reduced alternating diagram for a given knot K, then $D'$ is obtained from $D$ by a series of flypes by Menasco and Thistlethwaite's proof of the Tait Flyping conjecture in \cite{MT}. It is easy to check that a flype exchanges one crossing in a given diagram for another crossing of the same sign in the diagram resulting from the flype. This implies the following:

\begin{corollary}
The slopes of the checkerboard surfaces for any two reduced alternating diagrams for a given knot agree.
\end{corollary}

Thus, while our construction depends upon a specific choice of diagram for the alternating knot, the slopes of the checkerboard surfaces do not.

That these surfaces realize the maximal and minimal boundary slopes for a two-bridge or an alternating Montesinos knot now follows from the work of \cite{MMR} and \cite{IM1}. We show more generally that these surfaces realize the maximal and minimal integral boundary slopes for any alternating knot. The results for 2-bridge and alternating Montesinos knots are then immediate, since all boundary slopes for such knots are integral.
To prove that these surfaces realize the extreme integral slopes for arbitrary alternating knots, we must introduce the so-called {\em basic layered surfaces} of Adams and Kindred \cite{AK}.

Let $K$ be an alternating knot, and let $D$ be a reduced alternating diagram for $K$. Each crossing of $D$ has two possible smoothings, which we call {\em positive} and {\em negative} according to the convention in Figure 3.
The \textit{positive smoothing of the diagram $D$} is the smoothing in which each crossing is smoothed positively; the \textit{negative smoothing of the diagram $D$} is the dual smoothing, in which each crossing is smoothed negatively.

\begin{figure}[h]
\begin{center}
\leavevmode\hbox{}
\includegraphics[scale=1]{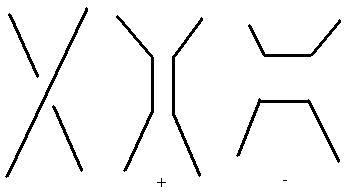}
\fcaption{{A crossing, its positive smoothing, and its negative smoothing}} \label{smoothings}
\end{center}
\end{figure}

If we smooth all crossings in the diagram $D$, the knot is split into a collection of circles, known as the \textit{state circles} of the smoothing of the diagram. We say that the smoothing of the diagram is \textit{adequate} if for each crossing in $D$, the two segments replacing the crossing strands in the smoothing are in different state circles. Note that both the positive and negative smoothings of $D$ are adequate by Theorem 9.5.4 of \cite{C} since $K$ is alternating.

\begin{definition}
A {\em basic layered surface} for an alternating knot $K$ is constructed by choosing any adequate smoothing of a reduced alternating diagram for $K$,  lifting nested state circles to be at distinct heights relative to the projection plane, filling each circle with a disk bounded in a plane parallel to the projection plane, and attaching half-twisted bands joining the disks at each crossing with twists chosen to agree with the crossing, as in Seifert's algorithm, so long as the bands can be attached without intersecting any of the disks.
\end{definition}
Note that the basic layered surfaces corresponding to the positive and negative smoothings of a reduced alternating diagram for $K$ are precisely the checkerboard surfaces for $K$ given by $D$.

As part of their main theorem, Adams and Kindred prove
\begin{theorem}(Adams and Kindred \cite{AK}, Theorem 5.3)
If $K$ is an alternating knot, then any spanning surface of $K$ has the same slope as one of the basic layered surfaces for $K$.
\end{theorem}

We obtain the following as an easy corollary of their work:
\begin{corollary}
If $K$ is an alternating knot with reduced alternating diagram $D$, then the two checkerboard surfaces for $K$ given by $D$ have the maximal and minimal boundary slopes among all essential surfaces with a single boundary component and an integral slope in the complement of $K$.
\end{corollary}

\begin{proof}
Let $K$ be any alternating knot, and let $S$ be an essential surface in the complement of $K$ with an integral slope and a single boundary component.  Then $S$ is a spanning surface for $K$, so by the theorem its slope agrees with that of one of the basic layered surfaces for $K$.  Adams and Kindred show in \cite{AK}, Proposition 3.2, that the slope of a basic layered surface is equal to $a - b + cr_+ - cr_-$, where $a$ is the number of crossings which were smoothed positively in the smoothing of $D$ and $b$ is the number of crossings which were smoothed negatively.  Clearly this is maximal and equal to $2 cr_+$ for the positive smoothing of $D$, for which $a$ is equal to the crossing number of $K$ and $b$ is zero, and it is minimal and equal to $2 cr_-$ for the negative smoothing of $D$, for which $a = 0$ and $b$ is the crossing number of $K$.  These smoothings yield the checkerboard surfaces for $K$ given by $D$.
\end{proof}

Finally we recall
\begin{theorem}(Ichihara and Mizushima \cite{IM1} and \cite{IM2})
If $K$ is an alternating Montesinos knot with reduced alternating diagram $D$, then the two checkerboard surfaces for $K$ given by $D$ have the maximal and minimal boundary slopes among all essential surfaces in the complement of $K$.
\end{theorem}


\section{Main Theorem}
Our final section is devoted to the proof of our main result:

\begin{theorem}
Let $K$ be an alternating knot. Then the largest boundary slope of an essential spanning surface of $K$ is twice the maximum degree of the Jones polynomial plus $\sigma(K)$ and the smallest boundary slope of an essential spanning surface of $K$ is twice the minimum degree of the Jones polynomial plus $\sigma (K)$, where $\sigma(K)$ denotes the signature of the knot $K$.
\end{theorem}

We begin by introducing known facts about the extreme degrees of the Jones polynomial.  Details may be found in Chapter 9 of \cite{C}.

Let $D$ be any diagram of any knot. Let $|S_\pm D|$ denote the number of state circles is the positive and negative smoothings of $D$, respectively.
Then the degree of any term of the Jones polynomial is bounded above and below by
$$ \frac{1}{2}(2cr_+ - cr_- + |S_-D| -1 ) \text{ and } \frac{1}{2}(cr_+ - 2cr_- - |S_+D| +1 ),$$
respectively. (This is an easy consequence of Theorem 9.5.1 of \cite{C}, which gives corresponding bounds for the extreme degrees of the bracket polynomial.) For any diagram for which the positive and negative smoothings of the knot are adequate, including any reduced alternating diagram of an alternating knot, these values are attained by terms in the Jones polynomial and therefore represent the highest and lowest degrees of the Jones polynomial of the knot.

Next, we turn to knot signatures. We review Gordon's and Litherland's characterization of the knot signature in terms of the signature of the Goeritz matrix \cite{GL}.

Given a knot diagram $D$, checkerboard the diagram so that the black region is bounded in the plane. To each crossing $p$ in $D$, assign a value $\eta(p) = \pm 1$ according to the convention in Figure 4. Let $\mu(D) = \sum{\eta(p)}$, where the sum is over all exceptional crossings $p$.  Now \begin{figure}
  \includegraphics[scale=1]{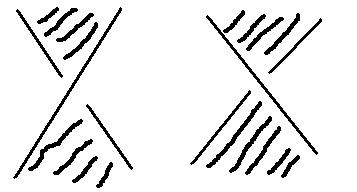}\\
  \caption{Define $\eta(p) = 1$ and $\eta(p)=-1$, respectively}\label{eta}
\end{figure}
label the white regions of the checkerboard $R_0, R_1, \dots, R_n$ so that $R_0$ is the unbounded white region.
\begin{definition}
The {\em Goeritz matrix} $G$ associated with $D$ is the $n \times n$ matrix given by $G = (G_{ij})$ where
$$
       G_{ij} = \left\{
         \begin{array}{lr}
           \sum_{p \in \partial R_i}{\eta(p)}  \text { if } i=j\\
           -\sum_{p \in \partial R_i \cap \partial R_j}{\eta(p)} \text{ otherwise}
         \end{array}
       \right.
$$
\end{definition}
Here, $i$ and $j$ range from 1 to $n$; $R_0$ does not contribute to either sum. Note that $G$ is symmetric, so it has a well-defined signature. Let $\sigma (G)$ denote the signature of $G$.

\begin{theorem} (Gordon and Litherland \cite{GL}, Theorem 6)
Let $D$ be a diagram of $K$ with associated Goeritz matrix $G$. Then
$$\sigma(K) = \sigma (G) - \mu(D) $$
where $\sigma(K)$ is the signature of the knot $K$.
\end{theorem}

We now prove Theorem 3.1. We remark that after we posted this paper to the arXiv, a closely related formula for the signature of an alternating knot was brought to our attention. This result of Traczyk was proved in the 1980's, appearing in \cite{T}.

\begin{proof}
Let $D$ be a reduced alternating diagram of $K$. Since $K$ is an alternating knot, by Corollary 2.9, the two checkerboard surfaces given by $D$ are the surfaces with the largest and smallest boundary slopes among all essential spanning surfaces of $K$. Hence by Proposition 2.5 the largest and smallest boundary slopes of essential spanning surfaces of $K$ are $2cr_+$ and $2cr_-$, respectively. Thus to prove our theorem it is sufficient to show
$$ cr_+ = |S_+D| - \sigma(K) -1 \text{ and } cr_- = |S_-D| + \sigma(K) -1. $$
We will prove the first equality; the second is proved analogously.

As in the proof of Proposition 2.5, if $S$ is a checkerboard surface given by $D$, then the exceptional crossings of $S$ are either all of the positive crossings or all of the negative crossings of $D$. Choose $S$ to be the surface for which the negative crossings are exceptional.

Replacing an outermost arc of $D$ with no crossings with an arc with the same endpoints wrapping around the outside of the diagram as necessary, we may assume that $S$ is built from regions which are bounded in the plane. We compute the signature of the Goeritz matrix $G$ for this diagram.

By our choices, $\eta(p) = -1$ for all crossings $p$.  Because $D$ is a reduced alternating diagram with $\eta (p) = -1$ for each crossing, it follows that that matrix $G$ is negative definite. To see this, one can easily check that for any $n$-vector $v$ we have $$v^T G v = - \sum_{i=1}^{n-1} \sum_{j=i+1}^n N_{ij} (v_i - v_j)^2 - \sum_{i=1}^n v_i^2 (N_i - \sum_{j \neq i} N_{ij}),$$
where $N_i$ is the number of crossings in $\partial R_i$ and $N_{ij}$ is the number of crossings in $\partial R_i \cap \partial R_j$. Since $N_{ij}\geq 0$  and $N_i - \sum_{j \neq i} N_{ij} \geq 0$ for all $i$ and $j$ and since $G$ is nonsingular, we see that $v^T G v <0$ for any vector $v$.

Thus, all eigenvalues of $G$ are negative, and the signature of $G$ is equal to the negative of one less then the number of white regions (since this is the size of the Goeritz matrix). By our choice of checkerboard, the number of white regions is just $|S_- D|$. Hence,
$$\sigma(G) = -|S_-D| +1. $$

Because $D$ is a reduced alternating diagram it is well know that $$cr(K) = |S_+D| +|S_-D| - 2.$$ (See \cite{C}, Lemma 9.4.3.) Then with $\mu(D)$ defined as above, we have
\begin{align*}
cr(K) + \mu(D) &= |S_+D| +|S_-D| +\mu(D) - 2 \\
&= |S_+D| - \sigma(G) +\mu(D) -1  \\
&= |S_+D| - \sigma(K) -1.
\end{align*}
Since $-\mu(D)$ is the number of exceptional crossings and a crossing is exceptional if and only if it is negative, $\mu(D) = -cr_-(D)$. Thus $ cr(K) + \mu(D)=cr_+$, which completes the proof.
\end {proof}

An interesting question for further study is to what extent boundary slopes are determined by the Jones polynomial for more general knots.  We note that for arbitrary knots it is not known that the checkerboard surfaces are essential or that they have extreme boundary slopes among spanning surfaces; there may be other essential spanning surfaces with greater or lesser slopes.  Also, for arbitrary knots the quantities
$$ \frac{1}{2}(2cr_+(D) - cr_-(D) + |S_-D| -1 ) \text{ and } \frac{1}{2}(cr_+(D) - 2cr_-(D) - |S_+D| +1 )$$
are bounds for the upper and lower degrees of the Jones polynomial which may not be realized. Thus the theorem proved here cannot be expected to hold in general.

\end{document}